\documentclass[aip,reprint,cha]{revtex4-1}   	
\usepackage[pdftex]{graphicx}

\usepackage[usenames,dvipsnames,svgnames,table]{xcolor}
\definecolor{black}{rgb}{0,0,0}
\definecolor{myblue}{rgb}{0,0,0.6}
\definecolor{mygreen}{rgb}{0,0.5,0}

\newcommand{\bo}{\mathbf}

\usepackage{amsmath,amssymb}
\usepackage{comment}
\usepackage{bbold}
\usepackage{bm}
\usepackage[normalem]{ulem}
\usepackage{tabularx}
\usepackage{multirow}
\usepackage{subfigure}

\DeclareMathOperator*{\argmin}{arg\,min}
\def\ul{\underline}

\newif\ifcolor
\colortrue

\begin{document}

\title{
Robust and optimal sparse regression for nonlinear PDE models
}
\author{Daniel R. Gurevich}
\affiliation{School of Physics, Georgia Institute of Technology, Atlanta, GA 30332, USA}
\author{Patrick A. K. Reinbold}
\affiliation{School of Physics, Georgia Institute of Technology, Atlanta, GA 30332, USA}
\author{Roman O. Grigoriev}
\affiliation{School of Physics, Georgia Institute of Technology, Atlanta, GA 30332, USA}
\date{21 July 2019}						

\begin{abstract}
This paper investigates how models of spatiotemporal dynamics in the form of nonlinear partial differential equations can be identified directly from noisy data using a combination of sparse regression and weak formulation. Using the 4th-order Kuramoto-Sivashinsky equation for illustration, we show how this approach can be optimized in the limits of low and high noise, achieving accuracy that is orders of magnitude better than what existing techniques allow. In particular, we derive the scaling relation between the accuracy of the model, the parameters of the weak formulation, and the properties of the data, such as its spatial and temporal resolution and the level of noise.
\end{abstract}

%% showpacs class option if PACS display desired
%% copied from siminos/blog/strategy.tex
% \PACS 02.20.-a \sep 05.45.-a \sep 05.45.Jn \sep 47.27.ed \sep 47.52.+j
% 02.20.-a      Group theory, mathematics
% 05.45.-a      Nonlinear dynamics and chaos
% 05.45.Jn      High-dimensional chaos


\maketitle

\begin{quotation}
In recent years, data-driven discovery of mathematical models of spatially extended systems described by nonlinear PDEs has emerged as a promising alternative to more traditional modeling approaches. Existing approaches to model discovery such as sparse regression have several major weaknesses, however. Most notably, they break down for data with high levels of noise and have to be tuned empirically to produce meaningful results, making them ill-suited for analyzing experimental data. We show how these weaknesses can be addressed using a weak formulation of the model PDE. The weak formulation has substantial freedom that makes it extremely powerful and flexible, but the question arises of how this freedom can be used to robustly obtain the most accurate model. This question is addressed here for the first time.
\end{quotation}

\section{Introduction}

Partial differential equations (PDEs) provide a natural description for the temporal evolution of spatially extended systems in various fields of science and engineering. 
Historically and practically important examples include wave equations arising in many areas of physics, the Schr\"{o}dinger equation in quantum mechanics, the Navier-Stokes equations in fluid dynamics, and reaction-diffusion equations used to model physical, chemical, or biological systems. 
In the past, models of such systems were almost always constructed from first principles or using a suitable empirical approach.
However, in recent years, a data-driven paradigm for learning the dynamics has emerged, which leverages the modern prevalence of data and computational power to create models when the underlying governing laws have eluded first-principles derivation.

Many indirect methods for learning the dynamics that do not require a PDE have been proposed.
Notable examples include equation-free modeling \cite{kevrekidis2003}, artificial neural networks \cite{hsu1997,raissi2018deep,pathak2018}, dynamic mode decomposition \cite{tu2013} and Koopman operator approaches \cite{mezic2013}, balanced truncation \cite{rowley2005}, and resolvent-based analysis \cite{mckeon2010}.
While these techniques can provide an economical approximate description of the dynamics, this is done at the cost of losing the mathematical structure that affords physical intuition or interpretability.
Symbolic regression, which was originally used to derive nonlinear ordinary differential equations describing low-dimensional systems \cite{bomgard2007, schmidt2009}, offers an enticing alternative by allowing construction of exact models and discovery of conservation laws.
The genetic algorithms used in these earlier studies are however computationally expensive, preventing application of this approach to high-dimensional systems.
Thus, the recent emergence of a sparse regression approach for model discovery \cite{xu_2008,khanmohamadi_2009,brunton2016,rudy2017,reinbold_2019} has made a significant impact.
Applied to spatially extended systems, this approach allows data-driven discovery of governing equations in the form of PDEs by evaluating a library of candidate terms containing partial derivatives at a large number of points and using a regularized regression procedure to compute the coefficients of each term and select a parsimonious model.

Sparse regression has proven computationally efficient and capable of reconstructing numerous canonical PDEs \cite{rudy2017,li_2019}, but it faces serious difficulties when used for analysis of experimental data.
One complication is that the proper choice of parsimonious model is often unclear. 
In many implementations, it relies on a manual Pareto analysis to balance model accuracy and complexity \cite{brunton2016} or on an automatic but complex thresholding procedure (e.g. sequential threshold ridge regression \cite{rudy2017}) that tends to be sensitive to the choice of parameters.
More importantly, existing sparse regression methods often suffer from low accuracy even in the absence of noise and completely break down at noise levels characteristic of realistic applications.
This is because they inherently require explicit numerical evaluation of partial derivatives of the data, which is a notably ill-conditioned problem.

In this paper, we present a weak formulation of the sparse regression problem that eliminates this fundamental issue.
We also suggest a simple thresholding procedure that can always identify the correct form of the governing PDE even in the presence of extremely high noise.
Finally, we explore how this extremely flexible and robust approach can be optimized and tuned to the properties of the underlying data set to maximize accuracy. 
This paper has the following structure.
Section \ref{sec:methods} describes our approach and the system used to test it.
Results are presented and interpreted in Section \ref{sec:results}, and conclusions are discussed in Section \ref{sec:conclusions}.

\section{Methods}
\label{sec:methods}

We consider the problem of using the data ${\bf u}({\bf x},t)$ to identify a parsimonious mathematical model in the form of a PDE
\begin{align}
\sum_{n=1}^N c_n \bo{f}_n(\bo{x},t,\bo{u}, \partial_t\bo{u},\partial^2_t\bo{u}, \nabla\bo{u},\nabla^2\bo{u}, \cdots) = 0
\label{eq:pde}
\end{align}
where each term in the sum is a function of ${\bf u}$ and its partial derivatives in space and time with constant coefficients $c_n$. 
In most applications, the form of the basis functions $\bo{f}_n$ can be restricted based on physical considerations, such as symmetries, conservation laws, etc. \cite{bar_1999,reinbold_2019}.
Typically, $\bo{f}_n$ are taken to be products of powers of independent variables (${\bf x}$, $t$) and dependent variables (${\bf u}$ and its various derivatives), although the form can be arbitrary in theory.
Our goal is to determine the constants $c_n$ for the terms that should be present in the model while eliminating the dynamically insignificant and thus likely spurious terms.
Sparse regression aims to convert the PDE (\ref{eq:pde}) to a tractable (and ideally, robust) linear algebra problem.
Conventionally this is done by evaluating all of the terms in the PDE at a random collection of points $({\bf x}_k,t_k)$ using finite differences \cite{vallette_1997,bar_1999}, spectral methods \cite{xu_2008,khanmohamadi_2009}, or polynomial approximation \cite{rudy2017,reinbold_2019}. 
All of these approaches are extremely sensitive to noise, especially when high-order derivatives are present.
We will instead pursue a weak formulation of the problem that can be obtained by multiplying (\ref{eq:pde}) by a weight $\bo{w}_j(\bo{x},t)$ and then integrating the result over a domain $\Omega_k$.
Repeating the process for $K$ distinct combinations of weight functions and integration domains yields the linear system
\begin{align}
Q\bo{c} = 0
\label{eq:lin}
\end{align}
where $\bo{c} = [c_1, \dots, c_N]^T$ and $Q = [\bo{q}_1, \dots, \bo{q}_N]$ is a ``library'' matrix, with each column $\bo{q}_n\in\mathbb{R}^K$ consisting of the integrals of the function $\bo{f}_n$ with all $K$ combinations of weights $\bo{w}_j$ and domains $\Omega_k$.

Note that there is an extra degree of freedom in \eqref{eq:lin} corresponding to the normalization of $\bo{c}$.
Conventionally this is dealt with by assuming that ${\bf f}_1=\partial_t{\bf u}$, setting $c_1=1$, and solving the overdetermined system that corresponds to the choice $K\gg N$ using least squares or some regularized version of it \cite{rudy2017}.
This is however not always a valid assumption: it is usually unknown {\it a priori} whether any given temporal derivative should be included in the PDE at all, whereas in this case a particular term is forced into the model.
Moreover, even if this term should be present in the model, the regression effectively assumes that the time derivative was computed without error, which reduces the practical accuracy of the procedure.

We will therefore not make the assumption that the model has the form of an evolution equation and consider the linear problem \eqref{eq:lin} in its most general form.
The normalization of $\bo{c}$ can be fixed by adding an extra row with arbitrary nonzero elements to $Q$, after which the resulting equation (\ref{eq:lin}) can be solved by 
ordinary least squares.
A more elegant solution pursued in the present study is to instead compute $\bo{c}$ 
as the right singular vector of $Q$ corresponding to the smallest singular value.
Note that this corresponds to the solution of a constrained least squares problem for $Q^TQ\bo{c}=0$:
\begin{align}
\bo{c}=\argmin_{\|\bo{c}\|=1}\|Q^TQ\bo{c}\|.
\label{eq:linn}
\end{align}
Once a suitable solution has been obtained by further constraining the problem, the resulting parsimonious model can be rewritten in the form of an evolution equation by solving for a term such as $\partial_t\bo{u}$ (or $\partial_t^2\bo{u}$ for a wave equation).

To obtain a parsimonious model, we employ an iterative procedure to eliminate unnecessary terms from \eqref{eq:pde}.
At each step $i$, singular value decomposition is used to obtain the solution ${\bf c}^i$ given the matrix $Q^i$, and the residual $\eta^i = \|Q^i {\bf c}^i\|$ is computed.
We then find the term with the smallest $\|c^i_n\bo{q}_n\|/\|\bo{q}_n\|$
and construct $Q^{i+1}$ by eliminating the column ${\bf q}_n$ from $Q^i$.
The corresponding term is eliminated from the model if $\eta^{i+1}<\gamma\eta^i$, where $\gamma>1$ is some fixed constant (we use $\gamma=1.4$ in the present study).
The iteration terminates at step $i$ if $\eta^{i+1}>\gamma\eta^i$, yielding a parsimonious model.
We find that this method compares favorably to alternatives such as sequential threshold ridge regression \cite{rudy2017} as it robustly eliminates spurious terms without requiring extremely careful choice of parameters.
Moreover, the sparsification parameter has a simple interpretation: $\gamma-1$ is the maximum acceptable relative increase in the residual resulting from discarding a single library term.

We illustrate the advantages of our approach by applying it to the Kuramoto-Sivashinsky equation \cite{kuramoto_1978, sivashinsky_1985}
\begin{align}\label{eq:kse}
c_1\partial_t u+c_2u\partial_x u+c_3\partial_x^2u+c_4\partial_x^4u = 0
\end{align}
which has posed a significant challenge in past studies of sparse regression \cite{xu_2008,rudy2017} because it contains a fourth-order partial derivative that is difficult to evaluate numerically with adequate accuracy.
Here $c_1=\cdots=c_4=1$ are all constants, although our approach can easily be extended even to the case when these coefficients are functions of time and/or space, as discussed below. 
Since this is a scalar equation in one spatial and one temporal dimension, we use scalar weight functions $w^j(x,t)$. 
If we denote the terms in the model \eqref{eq:kse} by $f_1,\cdots,f_4$, then
\begin{align}
&q_1^{jk} = \int_{\Omega_k} w_j\partial_t u\,d\Omega,\quad
 q_2^{jk} = \int_{\Omega_k} w_ju\partial_x u\,d\Omega, \nonumber\\
&q_3^{jk} = \int_{\Omega_k} w_j\partial^2_x u\,d\Omega,\quad
 q_4^{jk} = \int_{\Omega_k} w_j\partial^4_x u\, d\Omega,
\end{align}
where $d\Omega=dx\,dt$.
The key feature of the weak formulation is that it can almost always be used to completely eliminate, or at least reduce the order of, the derivatives acting on the noisy data by integrating by parts.
In our particular case, 
\begin{align}
&q_1^{jk} = -\int_{\Omega_k} u\partial_tw_j\, d\Omega,\quad 
 q_2^{jk} = -\frac{1}{2}\int_{\Omega_k} u^2\partial_xw_j\, d\Omega,\nonumber\\
&q_3^{jk} = \int_{\Omega_k} u\partial^2_x w_j\, d\Omega,\quad 
 q_4^{jk} = \int_{\Omega_k} u\partial^4_xw_j\, d\Omega
\end{align}
under the assumption that $w_j$ and its first three partial derivatives with respect to $x$ vanish on the boundary $\partial\Omega_k$.
In our implementation, we use the composite trapezoidal rule to evaluate the integrals numerically.

Note that although this particular PDE features constant coefficients, terms with variable coefficients can be treated in a similar manner.
For instance, suppose that the coefficient of the term $\partial_x^4u$ is a function of ${\bf x}$ and $t$ that can be expanded in some (finite) basis as
\begin{align}
c_4({\bf x},t) = \sum_p c'_pg_p({\bf x},t)
\end{align}
with some constants $c'_p$. Then
\begin{align}
\int_{\Omega_k} w_jc_4\partial^4_x u\, d\Omega = \sum_p c'_p q_p^{jk},
\end{align}
where 
\begin{align}
q_p^{jk}=\int_{\Omega_k} u\partial^4_x(g_pw_j)\, d\Omega.
\end{align}
Sparse regression for a model including such a term would then simply require expanding the library $Q$ to include additional columns ${\bf q}_p$ with entries $q_p^{jk}$.
In this case as well, no derivatives of the noisy $u$ are used in finding the elements of $Q$.

Although in principle integration domains of any shape can be used, here we will only consider rectangular domains of a fixed size
\begin{align}
\Omega_k = \{(x,t) \ : \ |x-x_k| \leq H_x, |t-t_k| \leq H_t\}
\end{align}
where the centers $(x_k, t_k)$ of the rectangles $\Omega_k$ are chosen randomly.
Similarly, there are many possible choices for the weight functions satisfying the boundary conditions on $\partial\Omega_k$; we focus on functions of the form
\begin{align}\label{eq:wf}
w_j = (\ul{x}^2-1)^\alpha(\ul{t}^2-1)^\beta e^{\pm il\pi\ul{x}}e^{\pm im\pi\ul{t}},
\end{align}
where $\ul{x} = (x-x_k)/H_x$, $\ul{t} = (t-t_k)/H_t$ are nondimensionalized independent variables and $\alpha\ge 4$, $\beta\ge 1$, $l\ge 0$, and $m \ge 0$ are integers.
Note that there are four weight functions (corresponding to the four different choices of the signs in the exponentials) for each pair of nonzero $l$ and $m$. The integrals $q_n^{jk}$ are all of the form
\begin{align}
F_n^{lm}=\int_{-1}^1d\ul{t}\int_{-1}^1d\ul{x} f_n^{\alpha\beta}(\ul{x},\ul{t}) e^{\pm il\pi\ul{x}}e^{\pm im\pi\ul{t}},
\end{align}
where
\begin{align}
f_n^{\alpha\beta}(\ul{x},\ul{t})=f_n (u,\ul{x},\ul{t}) (\ul{x}^2-1)^\alpha(\ul{t}^2-1)^\beta,
\end{align}
so $F_n^{lm}$ are the coefficients of the two-dimensional Fourier series for $f_n^{\alpha\beta}(\ul{x},\ul{t})$.
Although $f_n(u,\ul{x},\ul{t})$ is not periodic on $\Omega_k$, the functions $f_n^{\alpha\beta}(\ul{x},\ul{t})$ are.
Moreover, $f_n^{\alpha\beta}(\ul{x},\ul{t})$ has at least $\alpha-1$ continuous derivatives in $\ul{x}$ and ${\beta-1}$ continuous derivatives in $\ul{t}$, so the Fourier coefficients decay according to $F_n^{lm}\sim l^{-\alpha}m^{-\beta}$.
The powers $\alpha$ and $\beta$ therefore control the width of the Fourier spectrum of the entries $q_n^{jk}$ in the library $Q$, while the choice of $l$ and $m$ allows us to tune the frequencies of the weights to the spectral properties of the data.
The convergence rate of Fourier series turns out to control the accuracy with which the integrals are evaluated using data that are available only on a discrete grid.
For simplicity, we will assume that the same weight functions are integrated on every domain.
It is possible to use either weight functions involving only a single pair of frequencies (e.g., $l$ and $m$) or a range of frequencies in space and/or time.

To test our sparse regression approach, we computed a solution of the Kuramoto-Sivashinsky equation, using the integrator described in Ref. \onlinecite{rudy2017} to generate data on a physical domain with dimensions $L_x = 32\pi$ and $L_t = 500$.
The numerical integration generated data with spatial resolution $\Delta x=0.0491$ using a time step $\Delta t=0.005$, which was then downsampled to a lower spatial resolution $\delta_x$ and temporal resolution $\delta_t$.
Unless noted otherwise, the results presented below are for $\delta_x = 0.1964$ and $\delta_t = 1$.
For reference, the solution has a correlation length $\ell_x \sim1.67 \approx 8.5 \delta_x$ and correlation time $\ell_t\sim 8=8\delta_t$.
To test the effects of noise, Gaussian noise with standard deviation $\sigma s_u$ was added to the data for various choices of $\sigma$, where $s_u \approx 1.3$ is the sample standard deviation of $u$ on the whole domain.

\begin{figure*}[t]
\subfigure[]{\includegraphics[width=\columnwidth]{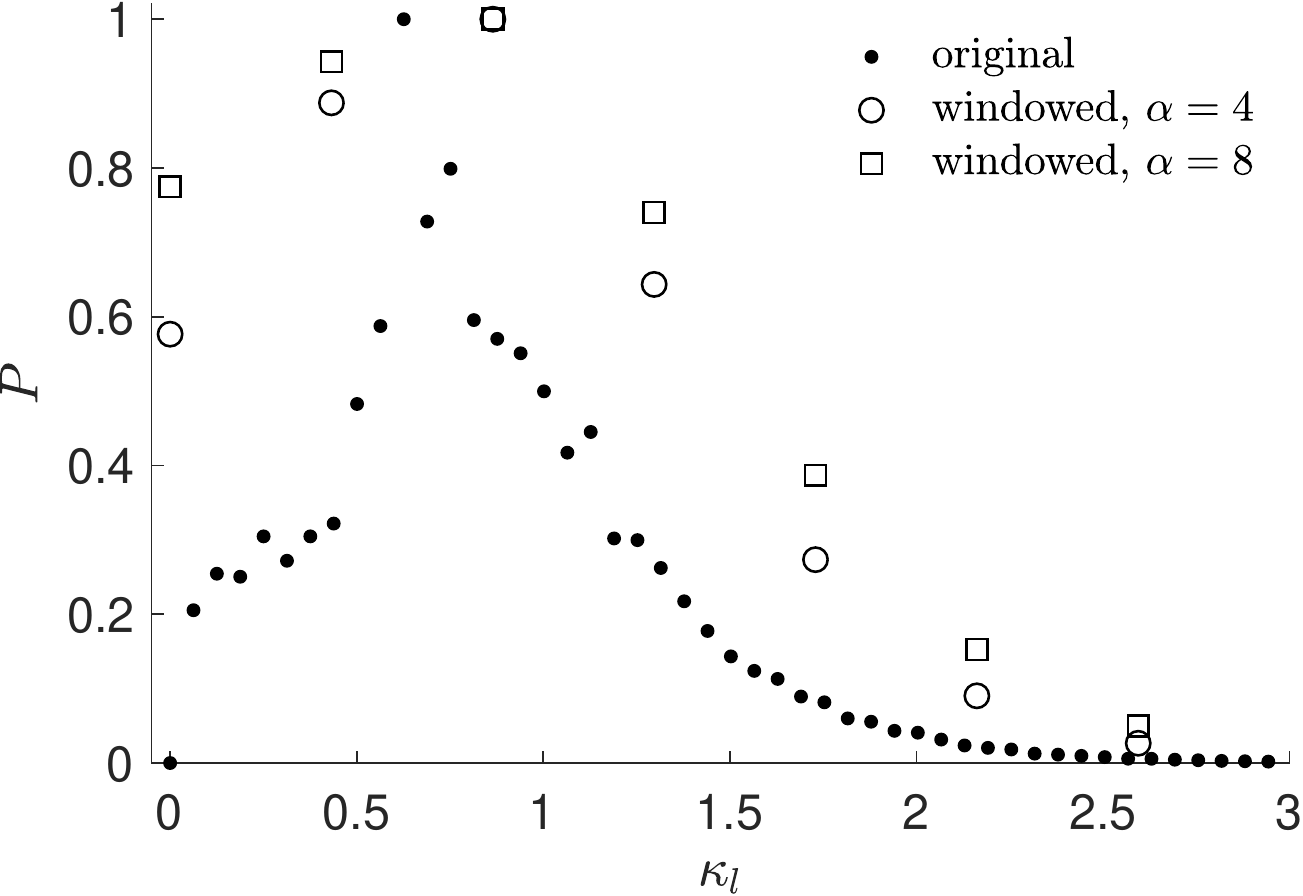}} \hspace{4mm}
\subfigure[]{\includegraphics[width=\columnwidth]{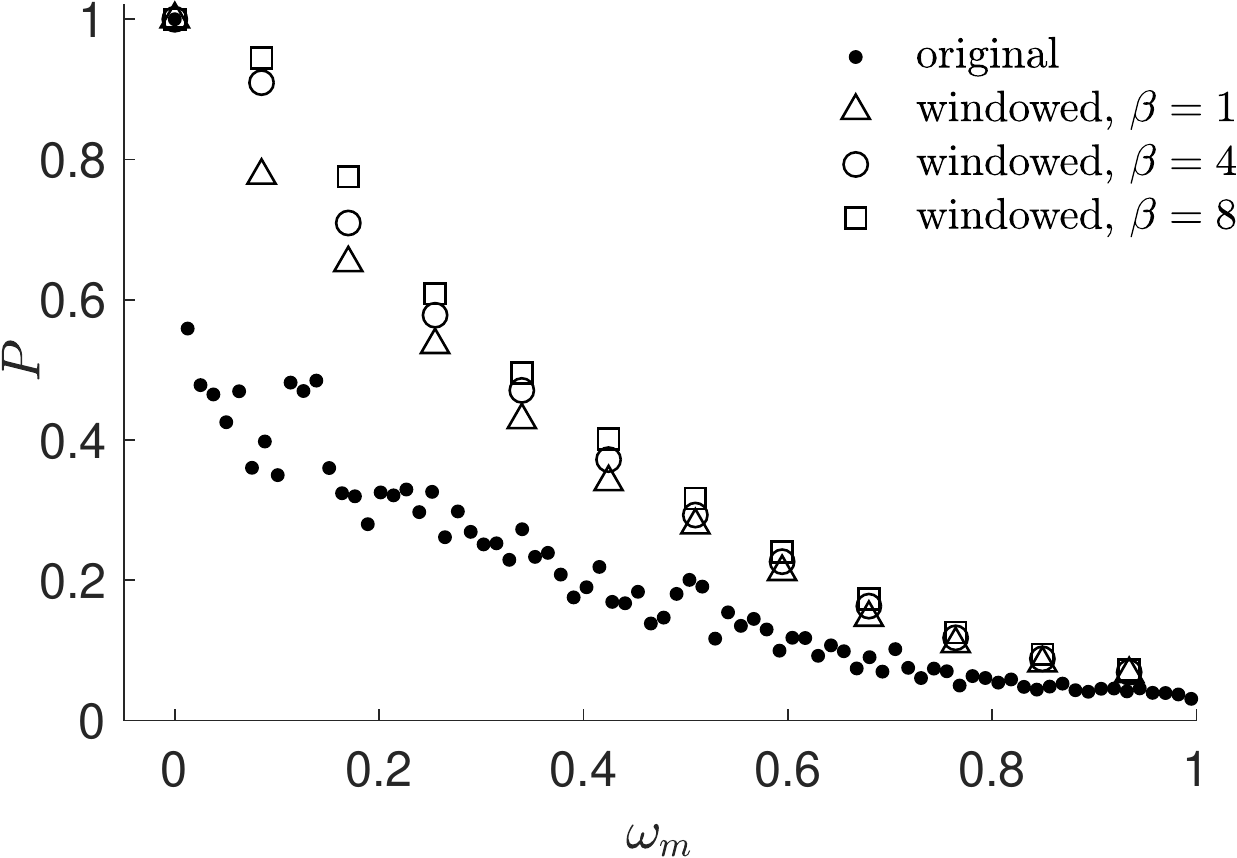}}
\caption{The power spectrum over (a) space and (b) time, normalized so that the maximum is 1. The black dots show the spectrum of the original data. The symbols correspond to the spectra of the windowed data multiplied by envelopes $E^{\alpha\beta}(\ul{x}, \ul{t})$ (with different choices of $\alpha$ or $\beta$) on a ``typical'' integration domain $\Omega_k$ (i.e., averaged over 1000 uniformly distributed choices of $(x_k, t_k)$. The spatial and temporal frequencies correspond to $\kappa_l=2\pi l/F_x$ and $\omega_m=2\pi m/F_t$ for the windowed data.
}
\label{fig:Fourier}
\end{figure*}

To test the ability of the algorithm to eliminate spurious terms, in addition to the terms present in the Kuramoto-Sivashinsky equation \eqref{eq:kse}, we also included terms $\partial_x u$, $\partial_x^3 u$, $u$, $u^2$, $u^3$, and $1$ (which represents a hypothetical forcing) in our library.
The corresponding integrals were rewritten using integration by parts to remove derivatives acting on $u$, as described previously. 
In the next section, we quantify the performance of our sparse regression approach using two key metrics: how well the algorithm can discriminate between the essential and spurious terms and how accurately it can determine the coefficients of the essential terms. 
Since the data were generated using a known model, we know which terms are essential (those contained in the PDE \eqref{eq:kse}).
If the reference model is unavailable, ensemble regression \cite{reinbold_2019} may be used instead to help distinguish essential terms from spurious ones.

\section{Results}
\label{sec:results}

As discussed previously, the elements of the library matrix $Q$ are given by the Fourier coefficients of the different terms included in the generalized model (windowed by the envelope $E^{\alpha\beta}(\ul{x}, \ul{t}) = (\ul{x}^2-1)^\alpha(\ul{t}^2-1)^\beta$ on each domain $\Omega_k$);
hence knowledge of the Fourier spectrum of the data is crucial for an optimal choice of the size of the integration domains $\Omega_k$ and the weight functions $w_j$.
The power spectrum (or, more precisely, the absolute value of the Fourier coefficients) of the noiseless data on the entire physical domain 
is shown in Figure \ref{fig:Fourier}.
In space, the spectrum is sharply peaked around a wave number $\kappa\approx 0.625$.
At high wave numbers, the spectrum decays exponentially, $P\propto e^{-\kappa/\bar{\kappa}}$ where $\bar{\kappa}\approx 0.3$.
In time, the spectrum is peaked at zero frequency $\omega$ and decays as a power law, $P\propto \omega^{-\chi}$ with $\chi\approx 2.5$.

Having characterized the data, we turn to the investigation of how the performance of our algorithm depends on the choice of various parameters.
Since the number of parameters is quite large, instead of exploring the entire parameter space, we focus on the dependence on one or two parameters at a time, with the remaining parameters staying fixed.
Specifically, the noise level $\sigma$ is fixed to 3\% and we use the following near-optimal parameters in the sparse regression. 
The dimensions of the integration domain are $F_x =2H_x=14.73$ and $F_t =2H_t=75$. This choice corresponds to an equal number of grid points in both directions, $F_x/\delta^*_x=F_t/\delta^*_t=75$.
Unless noted otherwise, we use a single set of weights with $\alpha=\beta=8$, $l=1$, $m=2$, and the sparsification parameter is $\gamma = 1.4$.
We generally use every combination of 4 weight functions over 50 integration domains, so that the total number of library rows is $K=200$.
To characterize the stochastic effects, for each set of parameters, we used an ensemble of $M=100$ trials featuring different random distributions of the integration domains and realizations of noise.

\begin{figure}[t]
\centerline{\includegraphics[width=\columnwidth]{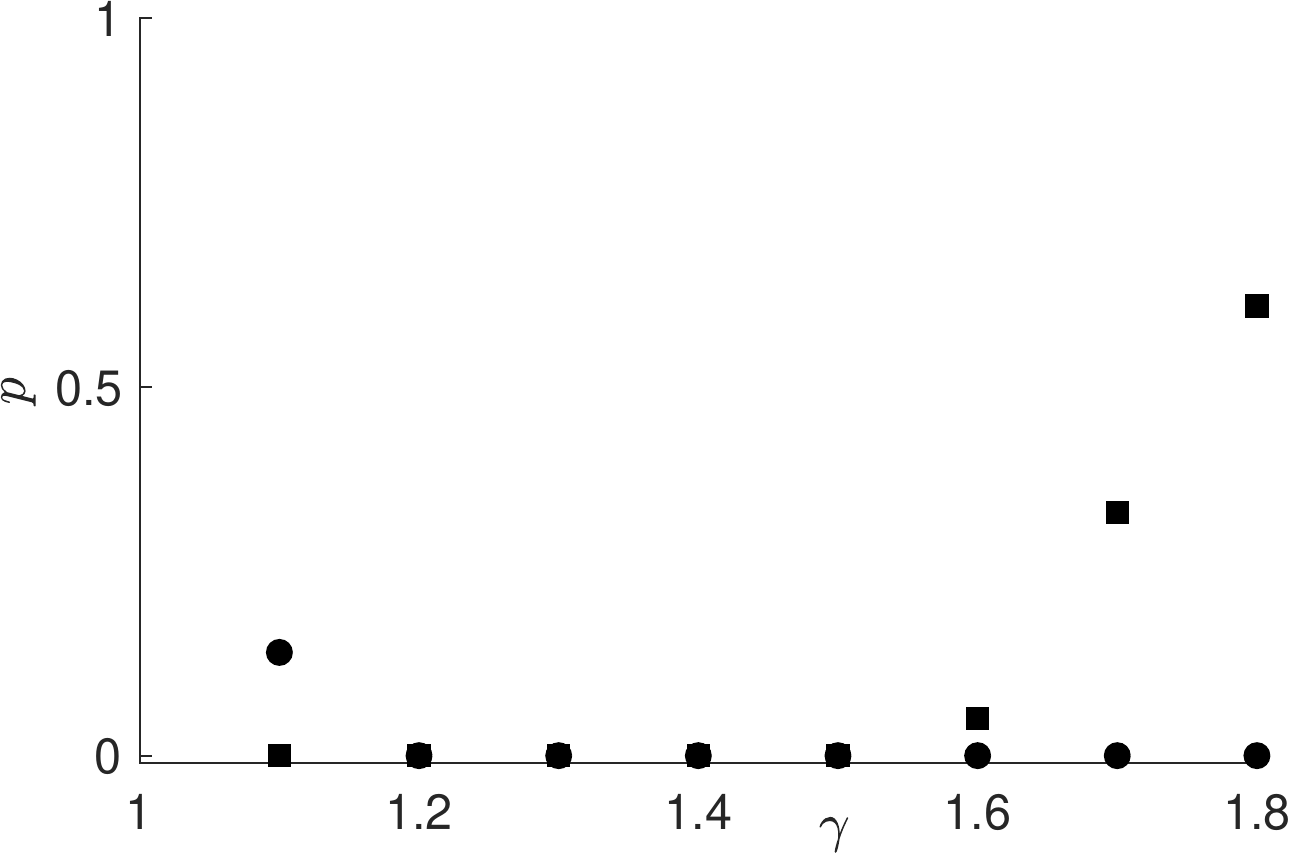}}
\caption{Fraction $p$ of identified models with spurious (circles) or missing (squares) terms at maximum noise level ($\sigma=1$) as a function of the sparsification parameter $\gamma$.
}
\label{fig:Spur}
\end{figure}

First, we tested the ability of the method to reconstruct the correct form of the PDE \eqref{eq:kse} for various values of $\gamma$ with all other parameters fixed at their near-optimal values.
Our iterative regression procedure proved very robust for a fairly wide range of values of $\gamma$. 
In particular, at a noise level of 30\%, it performed perfectly for $1.1\leq\gamma\leq2$, with the reconstructed model containing no missing or spurious terms in all of the trials.
For the highest noise level considered here (100\%), we found perfect performance for $1.2\leq\gamma\leq1.5$. In some fraction of the trials, spurious terms appeared at lower $\gamma$ and missing terms at higher $\gamma$, as shown in Figure \ref{fig:Spur}.
For reference, without the benefit of the weak formulation, sparse regression failed \cite{rudy2017} to correctly reconstruct the lambda-omega reaction-diffusion system, which is only second-order, for noise level as low as 1\%.

The accuracy of regression (i.e., model identification) was quantified by computing the relative error in each parameter of the Kuramoto-Sivashinsky equation
\begin{align}
\Delta c_n = \left|\frac{c_n-\bar{c}_n}{\bar{c}_n}\right|,
\end{align}
where $\bar{c}_n$ and $c_n$ are the true and estimated values of the model parameters, respectively, for $n=2,3,4$.
We normalize the estimated parameters so that $c_1 = 1$.
In all of the following figures, we plot the estimated mean value of $\Delta c_n$ with 95\% error bars, where all of the parameters in the regression procedure are held at their near-optimal values stated previously unless noted otherwise.

\begin{figure}[t]
\centerline{\includegraphics[width=\columnwidth]{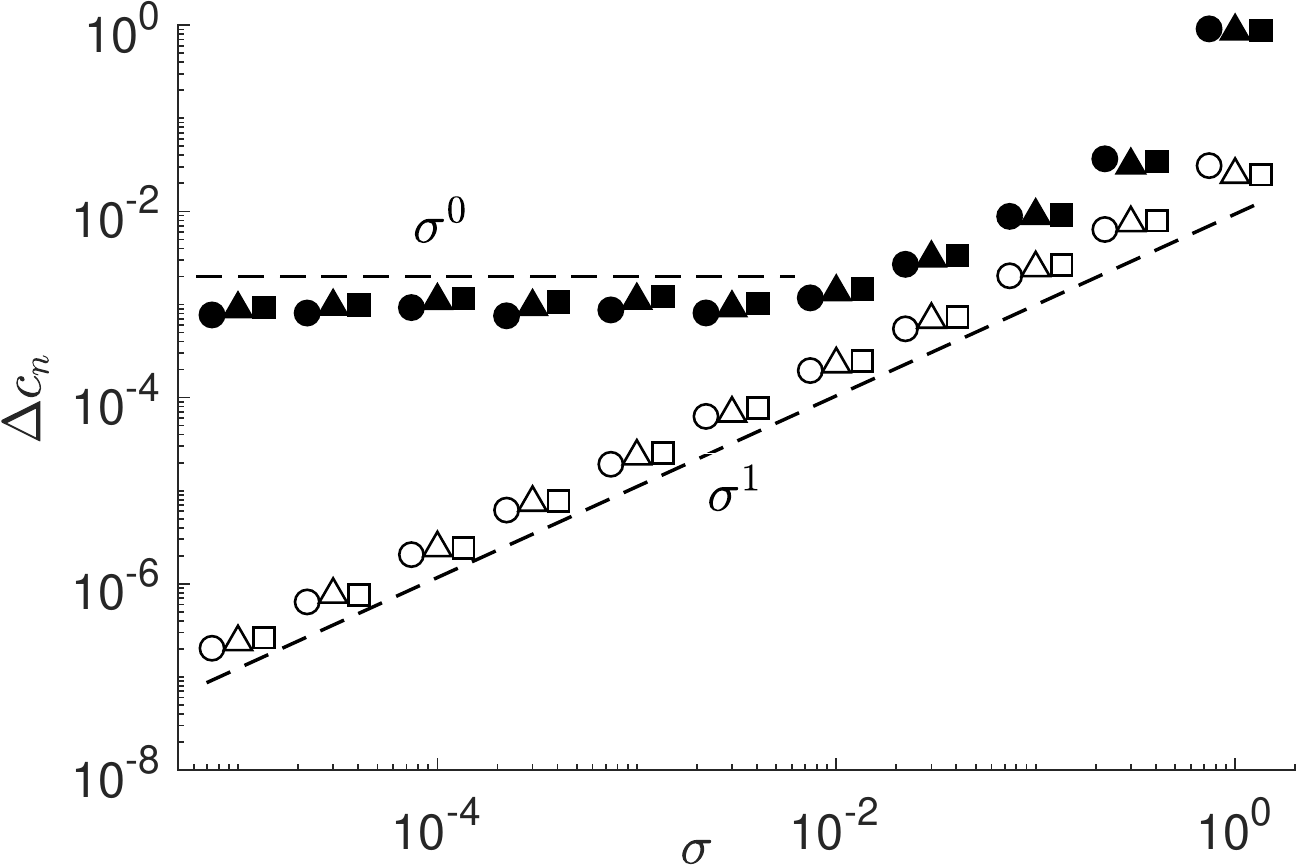}}
\caption{Parameter errors $\Delta c_n$ as a function of noise level. The circles, triangles, and squares correspond to $n=2, 3, 4$ respectively, and the empty and filled symbols indicate results for data with double resolution ($\delta_x = 0.0982, \delta_t = 0.5$) and half resolution ($\delta_x = 0.393, \delta_t = 2$), respectively. The dashed lines show the predicted scaling.
}
\label{fig:Noise}
\end{figure}

In particular, Fig. \ref{fig:Noise} shows the accuracy of regression as a function of $\sigma$ for two different choices of data resolution ($\delta_x$ and $\delta_t$).
It is worth noting that the average relative error is $\Delta c_n\sim 10^{-10}$ for all of the parameters for noiseless data with the higher of the two resolutions.
However, even for $1\%$ noise, $\Delta c_n\sim 2\times10^{-4}$, which is more than three orders of magnitude smaller than what had been achieved in previous studies \cite{rudy2017}.
The results are very similar for all three parameters; as this is generally the case, in subsequent figures, we only show the generally largest error $\Delta c_4$, which corresponds to the term $\partial_x^4 u$ involving the highest-order derivative.
We find two distinct regimes.
At higher noise levels, the error in evaluating the library matrix entries is due primarily to the averaged effect of noise. Applying the central limit theorem, we find that the relative error scales as
\begin{align}\label{eq:epsn}
\varepsilon_n\sim\sigma u_s\sqrt{\frac{\delta_x\delta_t}{F_xF_t}}.
\end{align}

At low noise levels, the parameter accuracy is controlled by numerical error, which has two different sources. 
The first source is a numerical error in the data itself, which is due to the finite accuracy of the integrator that ``solves'' the Kuramoto-Sivashinsky equation.
This source dominates for smaller $\delta_x$ and $\delta_t$.
For experimental data, this source would correspond to systematic error.
For larger $\delta_x$ and $\delta_t$, the parameter inaccuracy is mainly due to the error in computing the library matrix entries based on data that are available on a discrete grid.
Suppose we want to use numerical quadratures to evaluate an integral 
\begin{align}
I=\int_0^L g(x)dx,
\end{align}
where $g(x)\in C^{m}$ (i.e., has $m$ continuous derivatives) and $g^{(i)}(0)=g^{(i)}(L)$ for all $0\leq i<m$.
Then, for the composite trapezoidal rule on a grid with spacing $h$, the relative error associated with the discretization can be estimated using exact Euler-Maclaurin formulas\cite{trefethen2014} and is found to scale as $h^{m+2}|g^{(m+2)}|$ for $m$ even (or $h^{m+1}|g^{(m+1)}|$ for $m$ odd), where a characteristic value of the derivative on the interval $[0,L]$ is used.

Generalizing this result to two dimensions (and assuming $\delta_x\ll \min(\ell_x,F_x)$, $\delta_t\ll \min(\ell_t,F_t)$), we find an estimate of the relative discretization error for an element of the library matrix $Q$ that involves a temporal derivative of order $\nu_t$ and/or spatial derivative of order $\nu_x$:
\begin{align}\label{eq:mu}
\varepsilon_d\sim 
\begin{cases}
h^{\mu+2},& \mu\ \mathrm{even}\\
h^{\mu+1},& \mu\ \mathrm{odd}
\end{cases}
\end{align}
where $h=\delta_t/\ell_t\approx\delta_x/\ell_x$ and
$\mu=\min(\alpha-\nu_x,\beta-\nu_t)$.
It is easy to check that, due to the conditions on $\alpha$ and $\beta$, we always have $\mu\ge 0$, as it should be for the trapezoidal rule.
The Kuramoto-Sivashinsky equation features terms that all involve derivatives, with the lowest order being one and the highest being four;
hence, for even $\alpha=\beta\ge 4$, the exponent $\mu$ ranges between $\alpha-2$ and $\alpha+1$.
Therefore the scaling 
\begin{align} \label{eq:low_h}
\varepsilon_d\sim h^{\alpha-2}
\end{align} 
dominates for lower $h$, while the scaling
\begin{align} \label{eq:high_h}
\varepsilon_d\sim h^{\alpha}
\end{align}
dominates for higher $h$.

\begin{figure*}[t]
\subfigure[]{\includegraphics[width=\columnwidth]{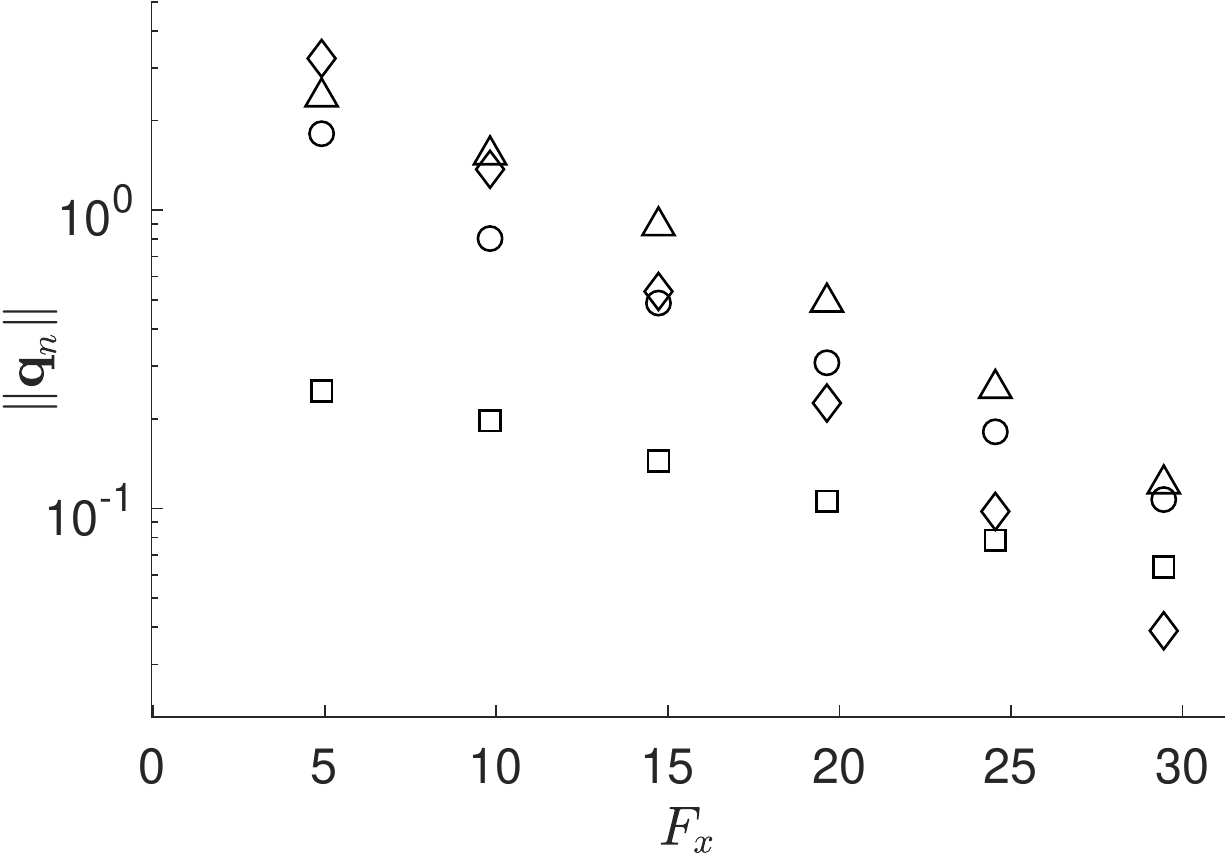}} \hspace{4mm}
\subfigure[]{\includegraphics[width=\columnwidth]{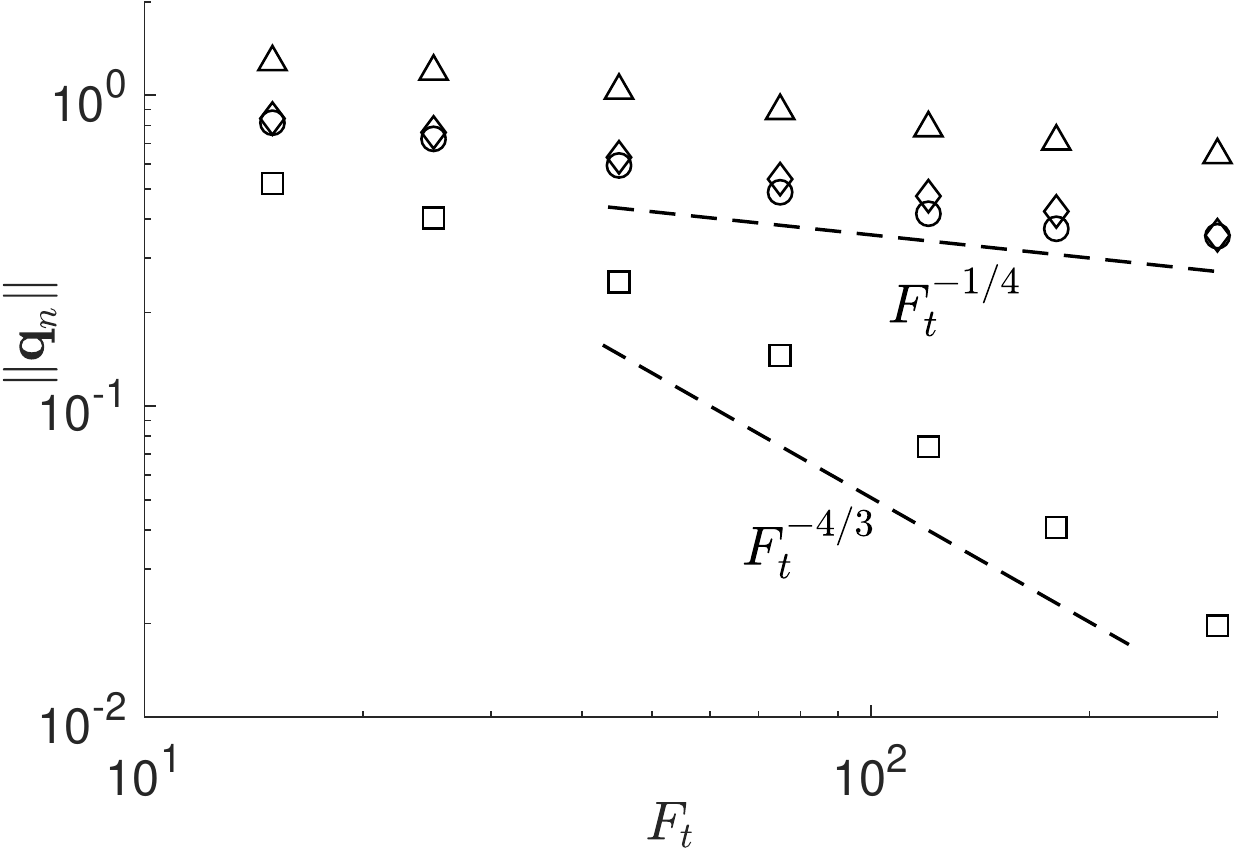}}
\caption{Scaling of the first four columns of the library $Q$ with the size of the integration domain in the (a) spatial and (b) temporal directions. The columns correspond to $\partial_t u$ (squares), $u\partial_x u$ (circles), $\partial_x^2 u$ (triangles), and $\partial_x^4 u$ (diamonds). 
}
\label{fig:norm}
\end{figure*}

\begin{figure}[t]
\subfigure[]{\includegraphics[width=\columnwidth]{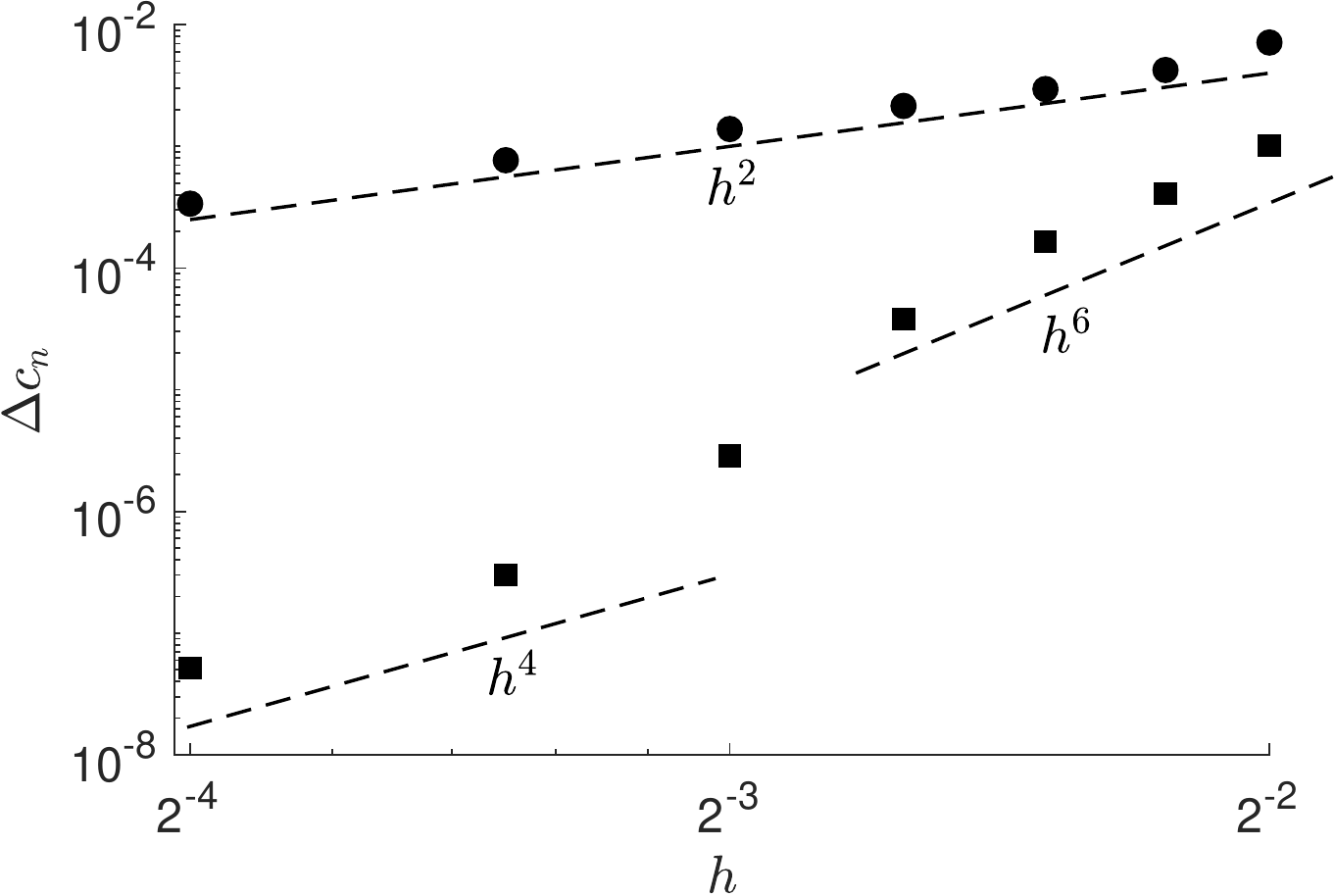}} 
\caption{Parameter error $\Delta c_4$ as a function of the resolution of noiseless data for $\alpha=\beta=4$ (circles) and $\alpha=\beta=6$ (squares). The dashed lines show the predicted scaling.
}
\label{fig:res}
\end{figure}

The error $\Delta c_n$ can be found using perturbation theory. 
Let $\bar{Q}$ be the library matrix evaluated using a continuous noiseless solution so that $\bar{Q}\bar{\bo{c}}=0$ exactly (we assume that $\bar{Q}$ corresponds to the parsimonious model).
In the presence of measurement noise and/or discretization error, the error in evaluating each entry $q_n^{jk}$ of the library matrix is proportional to $\varepsilon=\max(\varepsilon_d,\varepsilon_n)$, so
\begin{align}
Q=\bar{Q}+\varepsilon\hat{Q}
\end{align}
for some matrix $\hat{Q}$ whose entries are distributed as white Gaussian noise. 
Note that the entries of $\hat{Q}$ are $O(F_xF_t)$.
The entries of $\bar{Q}$ have a more complicated scaling that is determined by the Fourier spectrum of the data (i.e., exponential in space, power law in time). 
Specifically, we find (cf. Fig. \ref{fig:norm})
\begin{align}\label{eq:psin}
\|\bar{\bf q}_n\|\propto 
\begin{cases}
F_xF_t, & F_x\ll\ell_x,F_t\ll\ell_t\\
e^{-\lambda_nF_x}\left(\frac{\ell_t}{F_t}\right)^{\xi_n}\ell_x\ell_t,& F_x\gg\ell_x,F_t\gg\ell_t
\end{cases}
\end{align}
where $\lambda_n=O(\ell_x^{-1})$ and $\xi_n=O(1)$ are some positive constants.
To leading order in $\varepsilon$, the least squares solution to \eqref{eq:lin} is given by 
\begin{align}
\bo{c} = \bar{\bo{c}}-\varepsilon \bar{Q}^+\hat{Q}\bar{\bo{c}},
\end{align}
where $\bar{Q}^+$ is the Moore-Penrose pseudoinverse of $\bar{Q}$.
Since the elements of $\hat{Q}$ can be considered uncorrelated, we have for $F_x\gg\ell_x$ and $F_t\gg\ell_t$
\begin{align}\label{eq:pow0}
\Delta c_n \propto \varepsilon \frac{F_xF_tK^{-1/2}}{\psi(F_x,F_t)},
\end{align}
where the numerator and denominator describe the scaling of the entries of $\hat{Q}$ and $\bar{Q}$, respectively. Following from \eqref{eq:psin},
\begin{align}\label{eq:psi}
\psi(F_x,F_t)= e^{-\lambda F_x}\left(\frac{\ell_t}{F_t}\right)^{\xi}\ell_x\ell_t
\end{align}
with some positive constants $\lambda=O(\ell_x^{-1})$ and $\xi=O(1)$. 
For low $\sigma$, we have $\varepsilon=\varepsilon_d$ and therefore $\Delta c_n$ is independent of $\sigma$.
For high $\sigma$, we have $\varepsilon=\varepsilon_n$, so combining \eqref{eq:pow0} and \eqref{eq:epsn} we find
\begin{align}\label{eq:pow}
\Delta c_n \propto \sigma \sqrt{\frac{\delta_x\delta_t}{KF_xF_t}}\frac{F_xF_t}{\psi(F_x,F_t)}.
\end{align}
The predicted scaling of $\Delta c_n$ with $\sigma$ in both regimes is consistent with the results shown in Fig. \ref{fig:Noise}.
In particular, we find that the effect of changing the resolution of the data is quite minor at high $\sigma$, where $\Delta c_n\propto h$ according to \eqref{eq:pow}. At low $\sigma$, the effect is much stronger: for $\alpha=\beta=8$, we have $\Delta c_n\propto h^6$ according to \eqref{eq:mu}.
The dependence of the scaling in \eqref{eq:mu} on $\alpha$ and $\beta$ is further confirmed by Fig. \ref{fig:res}, which shows results for noiseless data.
In the $\alpha=\beta=4$ case, we observe the scaling law $\Delta c_n\propto h^2$ corresponding to \eqref{eq:low_h} in the entire range of $h$ we examined.
When $\alpha=\beta=6$, the parameter error scales according to $\Delta c_n\propto h^4$ for small $h$ and $\Delta c_n\propto h^6$ for large $h$, which correspond to the limiting cases \eqref{eq:low_h} and \eqref{eq:high_h}, respectively.
We should also note that for $h$ as large as $1/4$, the accuracy remains very good.
Thus, the method is suitable for fairly sparse data.

\begin{figure}[t]
\centerline{\includegraphics[width=\columnwidth]{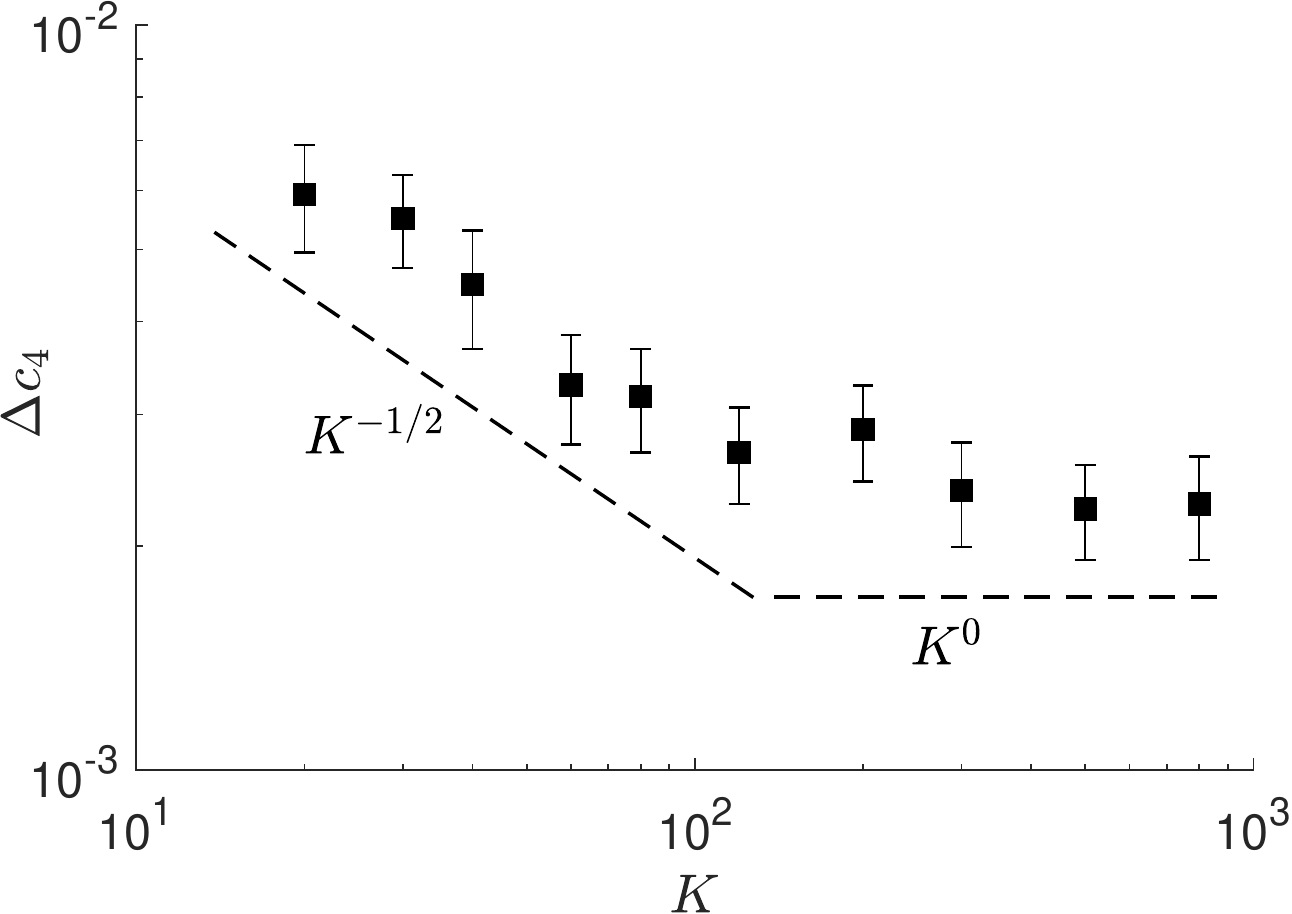}}
\caption{Parameter error $\Delta c_4$ as a function of the number of library rows $K$. Only the $l=m=0$ weight function is used and the physical domain size is reduced to $L_x=16\pi$, $L_t = 250$. The dashed lines show the predicted scaling.
}
\label{fig:numDom}
\end{figure}

As illustrated in Fig. \ref{fig:numDom}, we also observe the scaling for $\Delta c_n$ with $K$ predicted by \eqref{eq:pow}.
This scaling is expected to break down when the total area of the integration domains exceeds the area of the physical domain due to the loss of statistical independence between the data on different integration domains, leading to an increased linear dependence of the rows of the library matrix $Q$.
We can expect the error to asymptote to 
\begin{align}\label{eq:scaleNd}
\Delta c_n\propto\varepsilon N_d^{-1/2}
\end{align}
for $K\gg N_d$, where
\begin{align}
N_d=\frac{L_xL_t}{F_xF_t}
\end{align}
is the area ratio.
For the reference set of parameters, saturation did not occur over the range of $K$ we tested. 
To more easily observe the saturation effect, we set $l=m=0$, so that only one weight function is used and the number of integration domains equals $K$ (rather than $K/4$ for nonzero $l$ and $m$).
Furthermore, we reduce the size of the physical domain to $L_x=16\pi$ and $L_t = 250$, so that $N_d \approx 11$ is relatively small.
As Fig. \ref{fig:numDom} illustrates, for large $K$, the parameter accuracy indeed asymptotes to a constant.

\begin{figure}[t]
\includegraphics[width=\columnwidth]{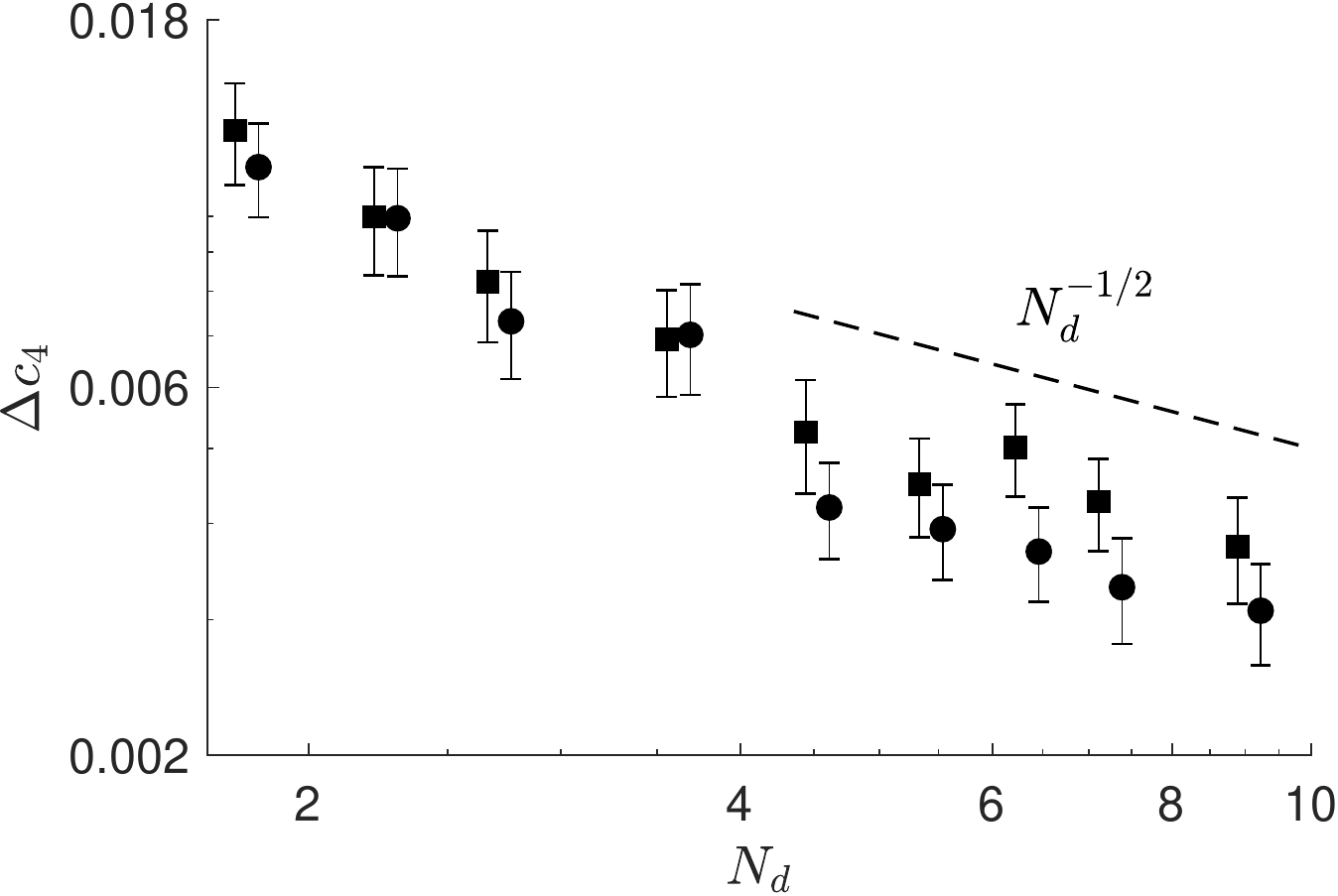}
\caption{Parameter error $\Delta c_4$ as a function of $N_d$ for $l=m=0$ and $K=500$. Squares correspond to fixing $L_t = 100$ and varying $L_x$ from $15.7$ to $98.2$. Circles correspond to fixing $L_x=19.6$ and varying $L_t$ from $80$ to $500$. The dashed line shows the predicted scaling.
}
\label{fig:L}
\end{figure}

The scaling described by \eqref{eq:scaleNd} can also be observed in the dependence of $\Delta c_n$ on the size of the physical domain (and hence $N_d$) with all other parameters fixed.
This dependence is quite important, since it determines how much data needs to be collected to identify the model with meaningful precision.
As Fig. \ref{fig:L} illustrates, choosing the physical domain to be just double the size of the (optimal) integration domain in both directions (which corresponds to $N_d=4$) already yields a rather acceptable accuracy when only one weight function is used.
When $l$ and $m$ are nonzero, accurate reconstruction is possible even if $N_d$ is only slightly greater than $1$.

\begin{figure*}[t]
\subfigure[]{\includegraphics[width=\columnwidth]{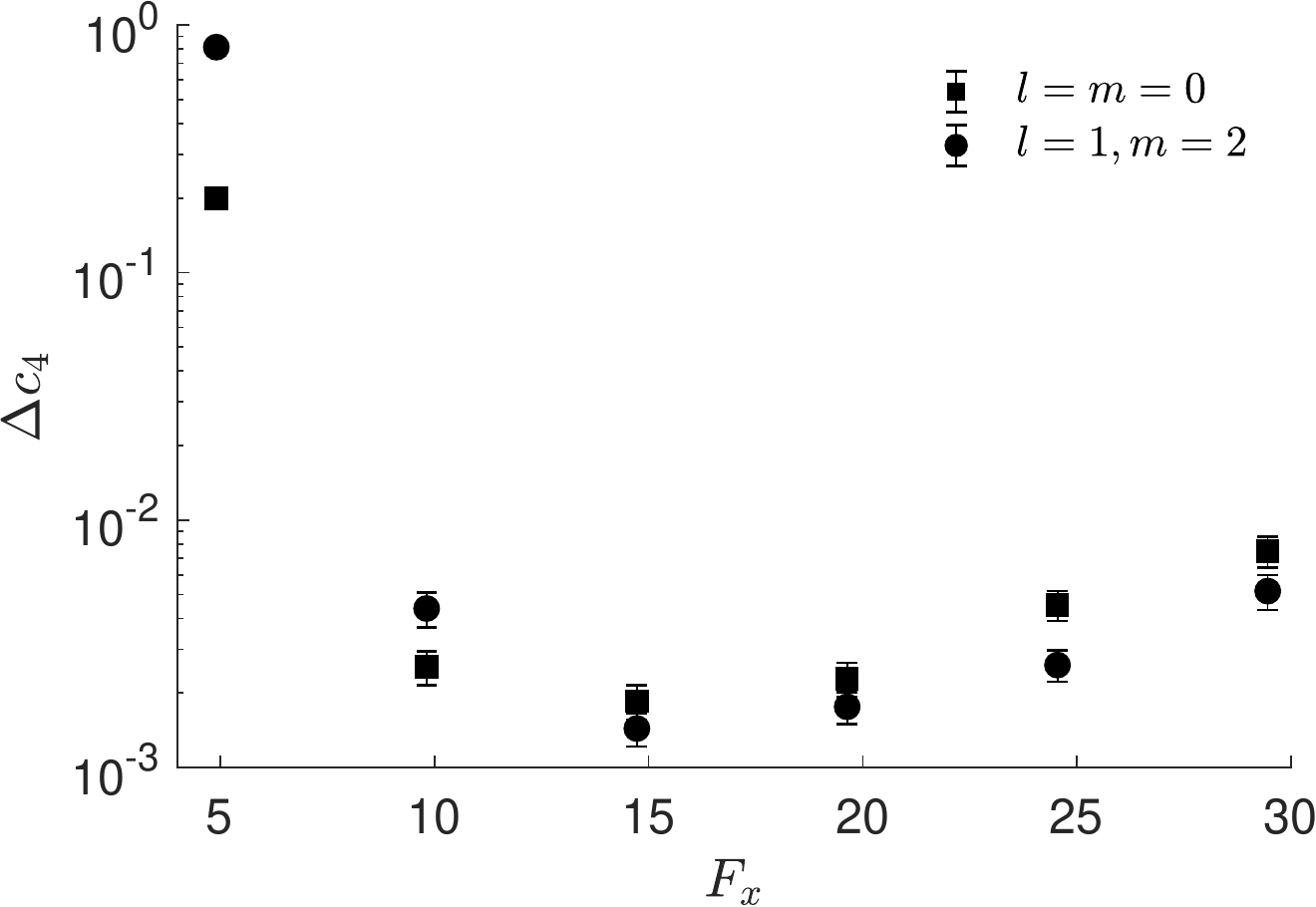}} 
\subfigure[]{\includegraphics[width=\columnwidth]{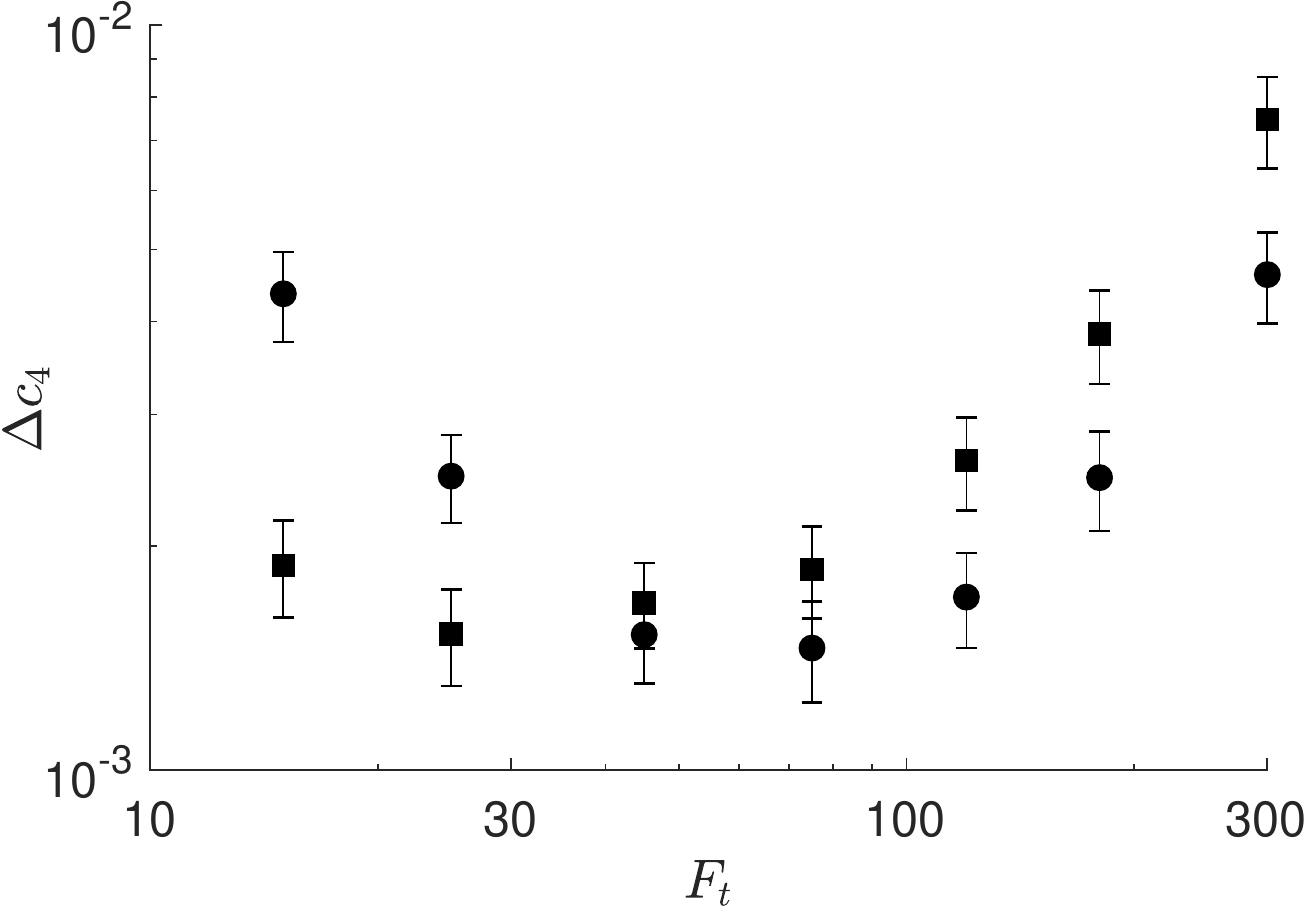}}
\caption{Parameter error $\Delta c_4$ as a function of the (a) spatial and (b) temporal dimensions of integration domains when only the $l=m=0$ weight function is used (squares) and for the optimal choice of $l$ and $m$ (circles).
}
\label{fig:size}
\end{figure*}

Next, we consider how the error in the estimated coefficients depends on the choice of the integration domain size.
Figure \ref{fig:size} shows the dependence of the error $\Delta c_4$ on the size of the integration domains for two different choices of the weight functions.
In panel (a), $F_t$ is fixed to 75 and $F_x$ is taken to vary, and in panel (b), $F_x$ is fixed to 14.73 with $F_t$ varying.
In both cases, we find that there is an optimal domain size with $F_x \approx 14.73$ and $F_t \approx 75$; moreover, the optimal values remain approximately the same even if we vary the size of the other dimension or the choice of weight functions.
For small $F_x$ and/or $F_t$, the error is large because (a) the integration domain is too small to effectively average out the influence of noise and (b) the numerical quadrature error is large (both $\epsilon_n$ and $\epsilon_d$ increase as $F_x$ and/or $F_t$ decrease). 
For large $F_x$ and $F_t$, we enter the regime described by \eqref{eq:pow0}, which predicts that the error should grow exponentially in $F_x$ and as a power of $F_t$.
Indeed, this is exactly what we observe in Fig. \ref{fig:size}.
Based on \eqref{eq:psin}, it appears that the optimal choice of $F_x$ and $F_t$ corresponds to the crossover between these two regimes, i.e., $F_x\propto\ell_x$ and $F_t\propto\ell_t$.
Our numerical results suggest that the optimal choice corresponds to $F_x/\ell_x\approx F_t/\ell_t\approx 8$.

\begin{figure*}[t]
\subfigure[]{\includegraphics[width=\columnwidth]{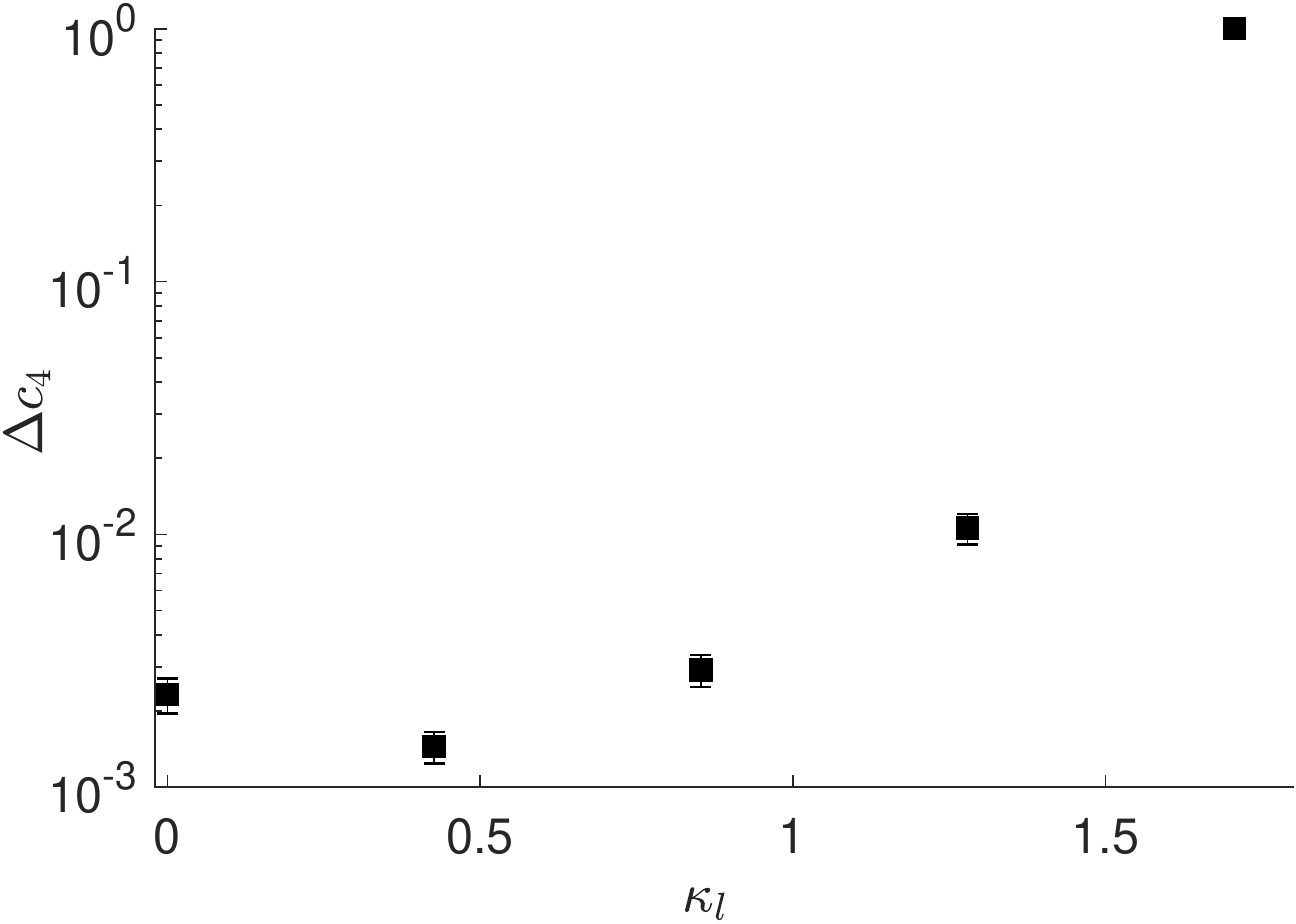}} 
\subfigure[]{\includegraphics[width=\columnwidth]{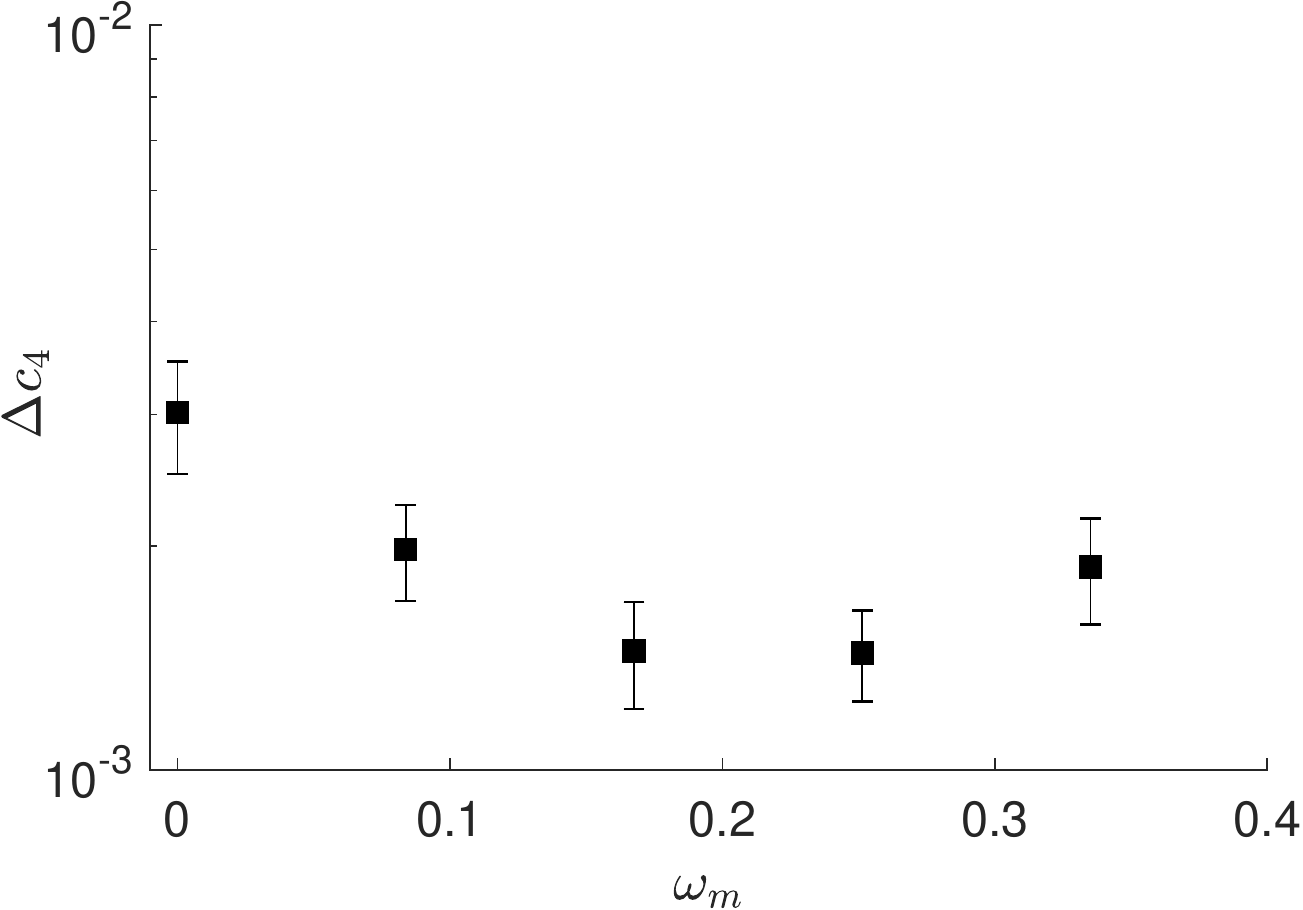}}
\caption{Parameter error $\Delta c_4$ as a function of (a) the wave number $\kappa_l= 2\pi l/F_x$ and (b) the frequency $\omega_m= 2\pi m/F_t$.
}
\label{fig:freq}
\end{figure*}

Finally, let us address the optimal choice of frequencies appearing in the weight functions \eqref{eq:wf}.
Figure \ref{fig:freq} shows the effect of varying either $l$ or $m$ with all other parameters fixed at their reference values.
Specifically, we plot $\Delta c_4$ versus $\kappa_l = 2\pi l/F_x$ and $\omega_m = 2\pi m/F_t$.
(Note that when $l$ or $m$ is 0, the number of distinct weight functions is halved, so we correspondingly double the number of integration domains to keep the number of rows in the library constant.)
One could assume that the optimal values would be given by the dominant frequencies of the original data (as we discussed previously, the dominant wave number is $\kappa\approx 0.625$ and the dominant temporal frequency is $\omega=0$).
According to Fig. \ref{fig:Fourier}, windowing the data broadens the peaks but leaves both dominant frequencies roughly the same: $\kappa \approx 0.8$ ($l=2$) and $\omega = 0$ ($m=0$).
Unfortunately, it turns out that we cannot use the spectra to exactly predict the optimal frequencies, which are $\kappa \approx 0.4$ ($l=1$) and $\omega \approx 0.2$ ($m = 2$ or $3$).
However, choosing the frequencies based on the spectra still produces reasonably good accuracy (within a factor of $4$ or so of the optimal result).

These results suggest that using weight functions with a combination of different frequencies may be more robust and/or accurate.
To test this hypothesis, we considered the case in which weight functions with a range of frequencies in space or time were included, with the total number of library rows fixed at $200$.
However, this approach yielded a decrease in the accuracy, as the broader choice of weight functions did not compensate for a decrease in the number of integration domains.
This suggests that the optimal strategy is to use a large number of integration domains while keeping the frequencies of the weight functions fixed.

\section{Conclusions}
\label{sec:conclusions}

We have introduced a robust and flexible approach to data-driven discovery of models in the form of nonlinear PDEs. 
The approach uses a weak formulation, coupled with a novel sparse regression procedure, to obtain a parsimonious description.
We have demonstrated its capability to identify PDEs, even with high-order derivatives, from extremely noisy data with unprecedented accuracy.
For instance, with 1\% noise, we were able to reduce the error in estimating the parameters of the 4th-order Kuramoto-Sivashinsky equation from 50\% \cite{rudy2017} to just $2 \times 10^{-4}$.
Furthermore, whereas correct identification of the functional form of the underlying PDE has been far from guaranteed at any noise level using past approaches, our algorithm was able to reconstruct the Kuramoto-Sivashinsky equation accurately in $100\%$ of cases from data with a signal-to-noise ratio of 100\%.

This impressive performance is achieved by shifting the partial derivatives from the data onto a known smooth weight function using integration by parts, thus avoiding the large errors incurred by repeated numerical differentiation.
Our method also proved to be well-adapted to sparse data, maintaining errors of less than $0.1\%$ for a grid resolution only $4$ times finer than the correlation length/time.
Such reliability and high accuracy in the presence of noisy or sparse data is indispensable for analysis of experimental data.
Notably, even in the absence of noise, our results compare very favorably with those of previous studies \cite{xu_2008,rudy2017} because the discretization error of the algorithm can be made extremely small: for the Kuramoto-Sivashinsky equation, the relative error in all parameters can easily be reduced to $10^{-10}$.
It is also important to mention that the computational cost of our algorithm is comparable to that of existing sparse regression methods.

We also derived the scaling laws that describe the accuracy of the regression as a function of the parameters used in the algorithm and the properties of the data.
These scaling laws can be used to fully exploit the flexibility of the weak formulation approach by tuning its various paramters.
In particular, the size of the input used by the regression can be controlled by choosing both the number of different integration domains and the number of different weight functions.
We have shown that the number of integration domains plays a much more important role than the number of weight functions: the best results can be obtained by using a set of weight functions with a fixed shape (frequency and envelope) and a large number of integration domains.
Furthermore, we have determined the optimal shape of the weights and the optimal size of the integration domains. 
The latter turned out to be determined by the correlation length and time describing the data (with the size roughly an order of magnitude larger than these characteristic scales).
We have also shown that, although the error can be reduced further by using data on ever-larger physical domains, satisfactory results can be obtained for physical domains that are just a factor of two larger than the optimal integration domain in each dimension.

\begin{acknowledgments}
This material is based upon work supported by the National Science Foundation under Grant No. CMMI-1725587. DG gratefully acknowledges the support of the Letson Undergraduate Research Scholarship.
\end{acknowledgments}

\section{References}
%\bibliography{ml}
 % Use this when bib file is not accepted

\end{document}